\documentclass[12pt,reqno]{amsart}
\usepackage{amsmath,amssymb,amsfonts,amscd,latexsym,amsthm,mathrsfs}
\usepackage{color}
\usepackage{relsize}
\usepackage{hyperref}
\usepackage{fullpage}
\usepackage{enumitem}  

\usepackage{mathtools}
\usepackage{tikz}

\newtheorem{theorem}{Theorem}%[section]

\newtheorem{lemma}[theorem]{Lemma}
\newtheorem{corollary}[theorem]{Corollary}

\title{Explicit estimates for the sum $\sum_{k=0}^{n} k! {n\choose k}^2 (-1)^{k}$}

\author{Anne-Maria Ernvall-Hyt\"onen}
\address{Anne-Maria Ernvall-Hyt\"onen, Mathematics and Statistics, P.O. Box 68, 00014 University of Helsinki, Finland}
\email{anne-maria.ernvall-hytonen@helsinki.fi}

\author{Tapani Matala-aho}
\address{Tapani Matala-aho, Mathematics, Aalto University, P.O. Box 11100, FI-00076 Aalto, Finland}
\email{tapani.matala-aho@aalto.fi}

\begin{document}

\maketitle

\begin{abstract}
We are interested in finding an explicit estimate to the binomial sum 
$Q_n(x)=\sum_{k=0}^{n} k! {n\choose k}^2 (-x)^{k}$ at $x=1$ for $n=0,1,2,\ldots$. 
Despite of its own interest the polynomial $Q_n(x)$ is important as the denominator in the Pad\'e identity of the Euler's factorial series 
$E(x) = \sum_{k=0}^{\infty} k! x^k$ as well as its close connection to a classical Laguerre polynomial
$L_n(x) = \frac{1}{n!} e^x \left(\frac{d}{dx}\right)^n (e^{-x}x^n)$. 
Our main result is the explicit bound 
$$\left|L_n(1)-\sqrt{\frac{e}{\pi}}\cdot \frac{\cos (2\sqrt{n}-\frac{\pi}{4})}{n^{1/4}}
+\frac{17}{48}\sqrt{\frac{e}{\pi}}\frac{\sin(2\sqrt{n}-\frac{\pi}{4})}{n^{3/4}}\right|<\frac{0.51}{n}$$ 
for all $n=0,1,2,\ldots$, which replaces the Fej\'er's asymptotic formula from 1909.
As a corollary of this, one also gets a new proof for the bound $|Q_{n}(1)| \le n!$, and even more. 
\end{abstract}

\section{Introduction}

We are interested in estimating the binomial sum 
\begin{equation}\label{summa}
\sum_{k=0}^{n} k! {n\choose k}^2 (-1)^{k}
\end{equation}
for $n=0,1,2,\ldots$.
The sum \eqref{summa} appears as a special  value of the Pad\'e polynomial 
\begin{equation}\label{padepol}
Q_n(x)= \sum_{k=0}^{n} k! {n\choose k}^2 (-x)^{k},
\end{equation}
the denominator polynomial in the Pad\'e identity $Q_{n}(x)E(x) - P_{n}(x)=R_{n}(x)$ of the Euler's factorial series 
$E(x) = \sum_{k=0}^{\infty} k! x^k$, see \cite{AMLOTA2023} and \cite{TAWA2018}. 

Our interest in the sum \eqref{summa} is motivated by the works \cite{AMLOTA2023} and \cite{TAWA2018}, which
consider arithmetic properties of the $p$-adic values the Euler's factorial series at integer points.
In particular, the sum \eqref{summa} is linked to the value $E(1) = \sum_{k=0}^{\infty} k!$, which is conjecturally
transcendental even its irrationality is an open problem. 
The Pad\'e approximation method used in \cite{AMLOTA2023} and \cite{TAWA2018} needs a good knowledge of the upper estimates
of the values of the polynomials $Q_n(x)$ and $P_n(x)$.
Let $x\in \mathbb{R} \setminus \{0\}$. From \cite{AMLOTA2023} we may find an explicit bound
\begin{equation}\label{QRAJA}
|Q_n(x)| \leq n!|x|^n e^{2\sqrt{\frac{n}{|x|}}},\quad n=0,1,2,\ldots.
\end{equation}
Although bound \eqref{QRAJA} is almost the best possible for $x<0$, it is not optimal for $x>0$, in particular at $x=1$.

A further motivation comes from the fact that the Pad\'e polynomial \eqref{padepol} is intimately related to a Laguerre polynomial
\begin{equation}\label{RODRIGUES}
L_n(x) := \frac{1}{n!} e^x \left(\frac{d}{dx}\right)^n (e^{-x}x^n) = \frac{1}{n!} \sum_{k=0}^{n} (n-k)! {n\choose k}^2 (-x)^{k}
= \sum_{k=0}^{n} \frac{1}{k!} {n\choose k} (-x)^{k}
\end{equation}
via 
\begin{equation}\label{EQUIVQL}
Q_n(x) = n!(-x)^n L_n(1/x). 
\end{equation}

Kooman and Tijdeman \cite{KoomanTijdeman1990} conjectured that $|L_{n}(1)| \le 1$ for all $n=0,1,2,\ldots$. 
Szeg\H{o} \cite{Szego1975}, formula (7.21.3), p.164, gives the bound $|L_n(1)|\leq e^{1/2}$. Furthermore, Korolev \cite{korolev}, 
Lemma 1.10, proved the bound $|L_n(1)|\leq 1$ while looking Gram's law at the zeros of the Riemann zeta function.

In fact, there is a more general result, Fej\'er's asymptotic formula,
\begin{equation}\label{AsympL}
L_{n}(x) = \frac{1}{\sqrt{\pi}n^{1/4}} \frac{e^{x/2}}{x^{1/4}} \cos\left(2\sqrt{nx}-\frac{\pi}{4}\right) 
+ \mathcal{O}\left(\frac{1}{n^{3/4}}\right),
\quad x>0, 
\end{equation}
proved already 1909 by Fej\'er \cite{Fejer1909}. See also Szeg\H{o} \cite{Szego1975}, Theorem 8.22.1, p.198, Tricomi \cite{Tricomi1949},
Kooman and Tijdeman \cite{KoomanTijdeman1990} and Borwein et al \cite{BBC2008}. 
In particular, \eqref{AsympL} shows 
\begin{equation}\label{AsympL1}
L_{n}(1) = \sqrt{\frac{e}{\pi}} \frac{1}{n^{1/4}} \cos\left(2\sqrt{n}-\frac{\pi}{4}\right) + \mathcal{O}\left(\frac{1}{n^{3/4}}\right).
\end{equation}
However, all these asymptotic results are non-explicit, meaning that we do not know the constant in the $\mathcal{O}$-symbol.
Of course, $L_{n}(1) \rightarrow 0$ as $n\rightarrow\infty$ by \eqref{AsympL1}.
But, for example, we can not deduce the fact, $|L_{n}(1)| \le 1$ for all $n=0,1,2,\ldots$, from \eqref{AsympL1}.   

First we will give an explicit version for \eqref{AsympL1}. 
\begin{theorem}\label{laguerre2} 
We have
\begin{equation}\label{paa2}
\left|L_n(1)-\sqrt{\frac{e}{\pi}}\cdot \frac{\cos (2\sqrt{n}-\frac{\pi}{4})}{n^{1/4}}+\frac{17}{48}\sqrt{\frac{e}{\pi}}\frac{\sin(2\sqrt{n}-\frac{\pi}{4})}{n^{3/4}}\right|<\frac{0.51}{n}
\end{equation}
for all $n\geq 1$.
\end{theorem}
To our knowledge, estimate \eqref{paa2} is the first explicit version of the asymptotic result \eqref{AsympL1}.

As an immediate corollary with the help of \eqref{EQUIVQL}, we obtain an equivalent result.
\begin{theorem}\label{laguerre2kertoma}
We have
\begin{equation*}\label{paa22}
\left|Q_n(1) - (-1)^n\left(\sqrt{\frac{e}{\pi}}\cdot \frac{n!}{n^{1/4}}\cos\left(2\sqrt{n}
-\frac{\pi}{4}\right)-\frac{17}{48}\sqrt{\frac{e}{\pi}}\frac{n!}{n^{3/4}}\sin\left(2r-\frac{\pi}{4}\right)\right)\right|<\frac{0.51n!}{n}
\end{equation*}
for all $n\geq 1$.
\end{theorem}

Finally, as a corollary of the previous results, we get a new proof for the bound $|L_n(1)|\leq 1$.
\begin{corollary}\label{laguerre3} 
We have
\begin{equation*}
|L_n(1)|\leq 1,
\end{equation*}
which is equivalent to
\begin{equation*}\label{paa1}
|Q_n(1)| \leq n!
\end{equation*}
for all $n\geq 0$.
\end{corollary}

\section*{Acknowledgements}

We would like to thank professors Juan Arias de Reyna, Maxim Korolev and Paul Voutier for valuable comments and references. 
Ernvall-Hyt\"onen's research is supported by the Emil Aaltonen foundation.

\section*{Structure of the proof}

The claims in our Theorems up to $n=10^5$ are verified by computations which are done by SAGE \cite{sage}. 
For that we note that the sequence $(Q_n(x))$ satisfies the recurrence
\begin{equation*}\label{POLRECURRENCES}
Q_{n+2}(x)  = (1-(2n+3)x) Q_{n+1}(x) - (n+1)^2x^2 Q_{n}(x),\quad n=0,1,\ldots,
\end{equation*}
with the initial values, $Q_0(x)=1$ and $Q_1(x)=1-x$. 
Therefore, we may use the recurrence
\begin{equation*}\label{sumrec-1}
Q_{n+2}(1)  = -2(n+1) Q_{n+1}(1) - (n+1)^2 Q_{n}(1),\quad n=0,1,\ldots,
\end{equation*}
with the initial values, $Q_0(1)=1$ and $Q_1(1)=0$, for fast computation of the values of $Q_{n}(1)$, $n=2,3\ldots$.

Computing
\[
n\left(\frac{(-1)^nQ_n(1)}{n!}-\sqrt{\frac{e}{\pi}}\cdot \frac{\cos (2\sqrt{n}-\frac{\pi}{4})}{n^{1/4}}+\frac{17}{48}\sqrt{\frac{e}{\pi}}\frac{\sin(2\sqrt{n}-\frac{\pi}{4})}{n^{3/4}}\right)
\]
by Sage and using the min and max commands, we notice that this difference is between $-0.1$ and $0.1$. Hence, the theorem is true for $n\leq 10^5$, and we may now concentrate on larger values of $n$.

For $n \ge 10^5$, we do the proof by analytic means.
Here the Rodrigues formula \eqref{RODRIGUES} implies a crucial contour integral representation
\begin{equation*}\label{}
L_n(x) =  e^x \frac{1}{2\pi i}\oint_{\gamma} \frac{e^{-z}z^n}{(z-x)^{n+1}}dz,
\end{equation*}
where $\gamma$ is a simple positively oriented loop (of winding number one) enveloping the point $z=x$.   

The proof of Theorem \ref{laguerre2} for $n\geq 10^5$, is inspired by ideas presented in Borwein et al \cite{BBC2008}.
 
The structure of the proof is the following:

\begin{enumerate}[label=(\roman*)]
    \item We wish to compute the value of 
    \begin{equation}\label{Ln1oint} 
    L_n(1)=\frac{e}{2\pi i}\oint_\gamma \frac{e^{-z}z^n}{(z-1)^{n+1}} dz,
    \end{equation}
    by choosing the path to be the circle with center at 
		$z=\frac{1}{2}$ and radius $r=\sqrt{n}$.
		We can limit the computations to large values of $n$ as small values can be checked with a computer.
    \item We change the variable and rewrite the integral \eqref{Ln1oint} in the form
\begin{equation*}
\begin{split}
L_n(1) & = \frac{\sqrt{e}r}{2\pi }\int_{-\pi}^{\pi} \frac{e^{-2ri\sin(\alpha)}}{-\frac{1}{2}+re^{i\alpha}}\cdot
			     \left(e^{-1/(re^{i\alpha})}\frac{1+\frac{1}{2re^{i\alpha}}}{1-\frac{1}{2re^{i\alpha}}}\right)^{n} e^{i\alpha}d\alpha \\
			 & = \frac{\sqrt{e}r}{2\pi }\int_{-\pi}^{\pi} \frac{e^{-2ri\sin(\alpha)}}{-\frac{1}{2}+re^{i\alpha}}\cdot
			     \left(e^{-t}\frac{1+t/2}{1-t/2}\right)^{n} e^{i\alpha} d\alpha,\quad t:=\frac{1}{re^{i\alpha}}.
\end{split}
\end{equation*}	
    \item To make the computation of the integral easier, we estimate the term 
		$\left(e^{-t}\frac{1+t/2}{1-t/2}\right)^{n}$ 
		with the polynomial $1+\frac{nt^3}{12}$.
    \item We are then left with the integral
$$\frac{\sqrt{e}r}{2\pi }\int_{-\pi}^{\pi} \frac{e^{-2ri\sin(\alpha)}}{-\frac{1}{2}+re^{i\alpha}}\cdot
\left(1+\frac{nt^3}{12}\right)
e^{i\alpha}d\alpha,$$ which we estimate by the integral $\frac{\sqrt{e}}{2\pi }\int_{-\pi}^{\pi} e^{-2ri\sin(\alpha)}d\alpha$. 
This estimate requires computing explicit coefficients and an explicit error term for the asymptotic expansion of the $J$-Bessel 
function using the approach by Watson \cite{watson}.
\item The integral $\frac{\sqrt{e}}{2\pi }\int_{-\pi}^{\pi} e^{-2ri\sin(\alpha)}d\alpha$ is the value $\sqrt{e}J_{0}(2r)$, 
which is then studied using an explicit expansion of the $J$-Bessel function.
\item Finally, all the error terms are put together and bounded.
\end{enumerate}

\section{Preliminaries}

\begin{lemma}\label{lemma4} 
Assume $|t|\le \frac{1}{100}$. Now
\begin{equation*}
\left|e^{-t}\frac{1+t/2}{1-t/2}-\left(1+\frac{1}{12}t^3+\frac{1}{80}t^5\right)\right| \le 0.00605|t|^6.			
\end{equation*}
\end{lemma}
\begin{proof}
We have
\begin{equation*}
\begin{split}
e^{-t}\frac{1+t/2}{1-t/2}
& = \left(1+\frac{t}{2}\right) \left(\sum_{k=0}^{\infty}\frac{(-t)^k}{k!}\right)  
    \left(\sum_{k=0}^{\infty}\left(\frac{t}{2}\right)^k\right) \\
& = \left(1+\frac{t}{2}\right) \sum_{k=0}^{\infty} a_k t^k = \sum_{k=0}^{\infty} b_k t^k
	= 1+\frac{1}{12}t^3+\frac{1}{80}t^5+\frac{1}{288}t^6+\frac{1}{448}t^7+\cdots,
\end{split}
\end{equation*}
where 
\begin{equation*}
\begin{split}
a_k & = \frac{1}{2^k}\left( \frac{2^0}{0!} - \frac{2^1}{1!} + \ldots +(-1)^k\frac{2^k}{k!} \right), \\
b_0 & = 1, \quad
b_k   = a_{k} + \frac{a_{k-1}}{2},\quad k=1,2,\ldots.
\end{split}
\end{equation*}
We also have
\begin{equation*}
a_k  = \frac{1}{2}a_{k-1} + \frac{(-1)^k}{k!},\quad 
b_k  = a_{k-1} + \frac{(-1)^k}{k!}.
\end{equation*}
Now we estimate
\begin{equation*}
\begin{split}
|a_k| & =   \frac{1}{2^k} \left| \sum_{m=0}^{k} (-1)^m\frac{2^m}{m!} \right| 
       \le \frac{1}{2^k} \left(\left| \sum_{m=0}^{\infty} (-1)^m\frac{2^m}{m!} \right| 
			                        + \left| \sum_{m=k+1}^{\infty} (-1)^m\frac{2^m}{m!} \right|\right) 
	     \le \frac{1}{2^k} \left( e^{-2} + \frac{2^{k+1}}{(k+1)!} \right). 
\end{split}
\end{equation*}
Hence
\begin{equation*}
\begin{split}
|b_k| & \le |a_{k}| + \left|\frac{a_{k-1}}{2}\right| 
        \le \frac{1}{2^k}\left( e^{-2} + \frac{2^{k+1}}{(k+1)!}\right) + \frac{1}{2^k}\left( e^{-2} + \frac{2^{k}}{k!}\right) \\
			& \le \frac{1}{2^k}\left( \frac{2}{e^2} + \left(1+\frac{2}{k+1} \right)\frac{2^{k}}{k!} \right) 
        \le \frac{1}{2^{k}}\left( \frac{2}{e^2} + \frac{9\cdot 2^{6}}{7\cdot 6!} \right) 
        \le \frac{0.38496}{2^k}
\end{split}
\end{equation*}
for $k\ge 6$. 
Finally
\begin{equation*}
\left| \sum_{k=6}^{\infty} b_k t^k \right| \le \frac{0.38496}{2^6}|t|^6 \sum_{k=0}^{\infty} \left|\frac{t}{2}\right|^k \le
\frac{0.38496}{2^6}\frac{200}{199}|t|^6<0.00605|t|^{6}.	
\end{equation*}
\end{proof}

\begin{lemma}\label{polynomiksi} Assume $|t|\leq \frac{1}{\sqrt{n}}$ and $|t|<\frac{1}{100}$. Now
\[
\left|\left(e^{-t}\frac{1+t/2}{1-t/2}\right)^{n}-\left(1+\frac{1}{12}t^3+\frac{1}{80}t^5\right)^{n}\right| 
<\frac{0.00605\cdot e}{n^2}.
\]
\end{lemma}

\begin{proof}
By Lemma \ref{lemma4} we know that
\[
\left|e^{-t}\frac{1+t/2}{1-t/2}-\left(1+\frac{1}{12}t^3+\frac{1}{80}t^5\right)\right| < 0.00605|t|^6.
\]
Hence
\begin{alignat*}{1}
&\left|\left(e^{-t}\frac{1+t/2}{1-t/2}\right)^{n}-\left(1+\frac{1}{12}t^3+\frac{1}{80}t^5\right)^{n}\right|\\
&=\left|e^{-t}\frac{1+t/2}{1-t/2}-\left(1+\frac{1}{12}t^3+\frac{1}{80}t^5\right)\right|\left|\sum_{k=0}^{n-1}\left(e^{-t}\frac{1+t/2}{1-t/2}\right)^k\left(1+\frac{1}{12}t^3+\frac{1}{80}t^5\right)^{n-1-k}\right|\\
& < 0.00605|t|^6\sum_{k=0}^{n-1}\left|\left(e^{-t}\frac{1+t/2}{1-t/2}\right)^k\left(1+\frac{1}{12}t^3+\frac{1}{80}t^5\right)^{n-1-k}\right|\\
& < 0.00605|t|^6n\left(1+\frac{1}{12}|t|^3+\frac{1}{80}|t|^5+0.00605|t|^6\right)^{n-1}\\ 
&\leq 0.00605|t|^4\left(1+\frac{1}{12}|t|^3+\frac{1}{80}|t|^5+0.00605|t|^6\right)^{n-1}
  < 0.00605\cdot e|t|^4\leq \frac{0.00605\cdot e}{n^2}
\end{alignat*}
since
\[
\frac{1}{12}|t|^3+\frac{1}{80}|t|^5+0.00605|t|^6<\frac{1}{n-1}.
\]
\end{proof}

\section{Proof of Theorem \ref{laguerre2} and Corollary \ref{laguerre3}}

Throughout the rest of the proof, assume $n\geq 10^5$.

We have
\[
L_n(1)=\frac{e}{2\pi i}\oint_{\gamma} \frac{e^{-z}z^n}{(z-1)^{n+1}} dz,
\]
where $\gamma$ is chosen to be the circle with center at $z=\frac{1}{2}$ and radius $r=\sqrt{n}$.

This proof is divided into three parts. In the first part, Preparation, we manipulate the integral slightly to turn it into form
\[
\frac{\sqrt{e}r}{2\pi }\int_{-\pi}^{\pi} \frac{e^{-2ri\sin(\alpha)}}{-\frac{1}{2}+re^{i\alpha}}\cdot\left(e^{1/(re^{i\alpha})}\frac{\frac{1}{2}+re^{i\alpha}}{-\frac{1}{2}+re^{i\alpha}}\right)^{n} e^{i\alpha}d\alpha.
\]
After this, in the second part, Reduction, we will first simplify the integral into a much simpler form
\[
\frac{\sqrt{e}r}{2\pi }\int_{-\pi}^{\pi} \frac{e^{-2ri\sin(\alpha)}}{-\frac{1}{2}+re^{i\alpha}}\cdot\left(1+\frac{nt^3}{12}\right) 
e^{i\alpha}d\alpha,\quad t:=\frac{1}{re^{i\alpha}}.
\]
by showing the the difference is sufficiently small. After this, in the last part we only have the main terms left, and we bound them.

\subsection{Changing the variable}

Write $z=\frac{1}{2}+\sqrt{n}\,e^{i\alpha}$, where $-\pi\leq \alpha\leq \pi $. Therefore
\begin{alignat*}{1}
L_n(1) = \frac{e}{2\pi i}\oint_{\gamma}\frac{e^{-z}z^n}{(z-1)^{n+1}} dz
       & =\frac{eir}{2\pi i}\int_{-\pi}^{\pi} \frac{e^{-(1/2+re^{i\alpha})}\left(\frac{1}{2}+re^{i\alpha}\right)^n}{\left(-\frac{1}{2}+re^{i\alpha}\right)^{n+1}}e^{i\alpha}d\alpha\\ 
			 &=\frac{\sqrt{e}r}{2\pi }\int_{-\pi}^{\pi} \frac{e^{-re^{i\alpha}}}{-\frac{1}{2}+re^{i\alpha}}\cdot\left(\frac{\frac{1}{2}+re^{i\alpha}}{-\frac{1}{2}+re^{i\alpha}}\right)^{n} e^{i\alpha}d\alpha.
\end{alignat*}
Now we may write
\begin{equation*}
\begin{split}
&\frac{\frac{1}{2}+re^{i\alpha}}{-\frac{1}{2}+re^{i\alpha}}=\frac{1+\frac{1}{2re^{i\alpha}}}{1-\frac{1}{2re^{i\alpha}}},\\
&e^{-re^{i\alpha}}=e^{-re^{i\alpha}-re^{-i\alpha}+re^{-i\alpha}}=e^{-re^{-i\alpha}-2ir\sin(\alpha)},\\
&t:=\frac{1}{re^{i\alpha}}.
\end{split}
\end{equation*}
Hence
\begin{equation*}
\begin{split}
L_n(1) & = \frac{\sqrt{e}r}{2\pi }\int_{-\pi}^{\pi} \frac{e^{-2ri\sin(\alpha)}}{-\frac{1}{2}+re^{i\alpha}}\cdot
			     \left(e^{-1/(re^{i\alpha})}\frac{1+\frac{1}{2re^{i\alpha}}}{1-\frac{1}{2re^{i\alpha}}}\right)^{n} e^{i\alpha}d\alpha \\
			 & = \frac{\sqrt{e}r}{2\pi }\int_{-\pi}^{\pi} \frac{e^{-2ri\sin(\alpha)}}{-\frac{1}{2}+re^{i\alpha}}\cdot
			     \left(e^{-t}\frac{1+t/2}{1-t/2}\right)^{n} e^{i\alpha} d\alpha.
\end{split}
\end{equation*}

\subsection{Approximating with a polynomial}

By Lemma \ref{polynomiksi}, we have
\[
\left|\left(e^{-t}\frac{1+t/2}{1-t/2}\right)^{n}-\left(1+\frac{1}{12}t^3+\frac{1}{80}t^5\right)^{n}\right| <\frac{0.00605\cdot e}{n^2}.
\]
Therefore we will approximate $L_n(1)$ by
\begin{equation*}
I_2:=\frac{\sqrt{e}r}{2\pi }\int_{-\pi}^{\pi} \frac{e^{-2ri\sin(\alpha)}}{-\frac{1}{2}+re^{i\alpha}}\cdot
\left(1+\frac{1}{12}t^3+\frac{1}{80}t^5\right)^{n} e^{i\alpha}d\alpha\,.
\end{equation*}
Now
\begin{equation}\label{E2}
\begin{split}
& |L_n(1)-I_2| \\ 
& = \frac{\sqrt{e}r}{2\pi }\left|\int_{-\pi}^{\pi} \frac{e^{-2ri\sin(\alpha)}}{-\frac{1}{2}+re^{i\alpha}}\cdot
\left(e^{-t}\frac{1+t/2}{1-t/2}\right)^{n} e^{i\alpha}d\alpha
- \int_{-\pi}^{\pi} \frac{e^{-2ri\sin(\alpha)}}{-\frac{1}{2}+re^{i\alpha}}\cdot
\left(1+\frac{1}{12}t^3+\frac{1}{80}t^5\right)^{n} e^{i\alpha}d\alpha\right| \\ 
& < \frac{\sqrt{e}r}{2\pi }\int_{-\pi}^{\pi} \left|\frac{e^{-2ri\sin(\alpha)}}{-\frac{1}{2}+re^{i\alpha}}\right|
\frac{0.00605\cdot e}{n^2} |e^{i\alpha}|d\alpha 
< \frac{0.00605\cdot e^{3/2}}{n^{3/2}(\sqrt{n}-1/2)} =: E_2.
\end{split}
\end{equation}
Let us continue simplifying the integral $I_2$. 
Look at the term $\left(1+\frac{1}{12}t^3+\frac{1}{80}t^5\right)^{n}$. We have
\[
\left|\left(1+\frac{1}{12}t^3+\frac{1}{80}t^5\right)^{n}-\left(1+\frac{1}{12}t^3\right)^{n}\right|=
\int_{1+\frac{1}{12}t^3}^{1+\frac{1}{12}t^3+\frac{1}{80}t^5}n|x|^{n-1}|\mathrm{d}x|\leq \frac{|t|^5n}{80}e\leq \frac{e}{80n^{3/2}}.
\]
So we introduce
\begin{equation*}
I_3:=\frac{\sqrt{e}r}{2\pi }\int_{-\pi}^{\pi} \frac{e^{-2ri\sin(\alpha)}}{-\frac{1}{2}+re^{i\alpha}}\cdot
\left(1+\frac{1}{12}t^3\right)^{n}
e^{i\alpha}d\alpha
\end{equation*}
and estimate
\begin{equation}\label{E3}
\begin{split}
 & |I_3-I_2| \\ 
      & \le \frac{\sqrt{e}r}{2\pi }\left|\int_{-\pi}^{\pi } \frac{e^{-2ri\sin(\alpha)}}{-\frac{1}{2}+re^{i\alpha}}
            \cdot\left(1+\frac{1}{12}t^3+\frac{1}{80}t^5\right)^{n} e^{i\alpha}d\alpha
					- \int_{-\pi}^{\pi } \frac{e^{-2ri\sin(\alpha)}}{-\frac{1}{2}+re^{i\alpha}}\cdot
				    \left(1+\frac{1}{12}t^3\right)^{n} e^{i\alpha}d\alpha\right| \\ 
      & <  \frac{\sqrt{e}r}{2\pi }\int_{-\pi}^{\pi} \left|\frac{e^{-2ri\sin(\alpha)}}{-\frac{1}{2}+re^{i\alpha}}
           \right|\frac{e}{80n^{3/2}} |e^{i\alpha}|d\alpha 
				< \frac{ e^{3/2}}{80n(\sqrt{n}-1/2)} =: E_3.
\end{split}
\end{equation}
Let us now approximate the term $\left(1+\frac{1}{12}t^3\right)^{n}$ by $1+\frac{nt^3}{12}$. We have
\begin{alignat*}{1}
\left|\left(1+\frac{1}{12}t^3\right)^{n}-1-\frac{nt^3}{12}\right|& \leq
\sum_{k=2}^{n}\binom{n}{k}\frac{|t|^{3k}}{12^k}\leq \sum_{k=2}^{n}\frac{n^k}{12^kn^{3k/2}}\\ 
&\leq \sum_{k=2}^{n}\frac{1}{12^kn^{k/2}}\leq \frac{1}{12^2n}\cdot \frac{1}{1-\frac{1}{12n^{1/2}}}<\frac{0.007}{n}.
\end{alignat*}
Thereby we write
\begin{equation*}
I_4:= \frac{\sqrt{e}r}{2\pi }\int_{-\pi}^{\pi} \frac{e^{-2ri\sin(\alpha)}}{-\frac{1}{2}+re^{i\alpha}}\cdot
\left(1+\frac{nt^3}{12}\right)
e^{i\alpha}d\alpha\,.
\end{equation*}
Then we estimate
\begin{equation}\label{E4}
\begin{split}
      & |I_4-I_3| \\ 
      & \le \left| \frac{\sqrt{e}r}{2\pi }\int_{-\pi}^{\pi } \frac{e^{-2ri\sin(\alpha)}}{-\frac{1}{2}+re^{i\alpha}} \cdot
                   \left(1+\frac{1}{12}t^3\right)^{n}
						        e^{i\alpha}d\alpha 
			           - \frac{\sqrt{e}r}{2\pi } \int_{-\pi}^{\pi } \frac{e^{-2ri\sin(\alpha)}}{-\frac{1}{2}+re^{i\alpha}}\cdot
				           \left(1+\frac{nt^3}{12}\right)
						        e^{i\alpha}d\alpha \right| \\ 						
      & <  \frac{\sqrt{e}r}{2\pi }\int_{-\pi}^{\pi } \left|\frac{e^{-2ri\sin(\alpha)}}{-\frac{1}{2}+re^{i\alpha}}\right|\cdot\frac{0.007}{n} 
           \cdot |e^{i\alpha}|d\alpha 
				<  \frac{0.007\sqrt{e}}{n^{1/2}(\sqrt{n}-1/2)} =: E_4.
\end{split}
\end{equation}

In the next chapter we will estimate $|I_4|$, which gives the main contribution.

\subsection{Treating the main terms}

In the following we will use the $J$-Bessel functions, which may be written as
$$
J_k(x) = \frac{1}{2\pi}\, \int_{-\pi}^{\pi} e^{-xi\sin(\alpha)+ki\alpha} d\alpha
$$
for $k\in\mathbb{Z}$. They satisfy the identity $J_{-k}(x) = (-1)^kJ_k(x)$ for $k\in\mathbb{Z}$,
see \cite{GradshteynRyzhik}, p. 912, formulae 8.411 and 8.404, respectively.

We wish to estimate the integral 
$I_4 = \frac{\sqrt{e}r}{2\pi }\int_{-\pi}^{\pi} \frac{e^{-2ri\sin(\alpha)}}{\frac{1}{2}+re^{i\alpha}}\cdot
\left(1+\frac{nt^3}{12}\right) e^{i\alpha}d\alpha$ by the integral 
$H_1 := \frac{\sqrt{e}}{2\pi }\int_{-\pi}^{\pi} e^{-2ri\sin(\alpha)}d\alpha=\sqrt{e}J_{0}(2r)$. 
Let us first estimate the error coming from this. Remember that $t=\frac{1}{re^{i\alpha}}$. Then look at the difference
\begin{equation*}\label{E56}
\begin{split}
I_4 - H_1 & = \frac{\sqrt{e}r}{2\pi }\int_{-\pi}^{\pi} \frac{e^{-2ri\sin(\alpha)}}{-\frac{1}{2}+re^{i\alpha}}\cdot
                \left(1+\frac{nt^3}{12}\right)
                e^{i\alpha}d\alpha -\frac{\sqrt{e}}{2\pi }\int_{-\pi}^{\pi} e^{-2ri\sin(\alpha)}d\alpha\\ 
            & = \frac{\sqrt{e}}{2\pi }\int_{-\pi}^{\pi} e^{-2ri\sin(\alpha)}\left(\frac{1}{1-\frac{t}{2}}\cdot
                \left(1+\frac{nt^3}{12}\right)-1\right)d\alpha \\
            & = \frac{\sqrt{e}}{2\pi }\int_{-\pi}^{\pi} e^{-2ri\sin(\alpha)}\left(
                \left(1+\frac{nt^3}{12}\right)\sum_{j=0}^{\infty}\left(\frac{t}{2}\right)^j-1\right)d\alpha\\ 
            & \leq  \frac{\sqrt{e}}{2\pi }\int_{-\pi}^{\pi} e^{-2ri\sin(\alpha)}\left(
                \left(1+\frac{nt^3}{12}\right)\sum_{j=2}^{\infty}\left(\frac{t}{2}\right)^j+\frac{nt^4}{24}\right)d\alpha\\
						& 	+ \frac{\sqrt{e}}{2\pi }\int_{-\pi}^{\pi} e^{-2ri\sin(\alpha)}
                   \left(\frac{nt^3}{12}+\frac{t}{2}\right)d\alpha
							=: E_5 + E_6.
\end{split}
\end{equation*}

Let us start with looking at the term $E_5$. We will simply use absolute values to bound the term. We will do it in two parts. 
The leftmost one will contribute
\[
\frac{\sqrt{e}}{2\pi }\int_{-\pi}^{\pi} \left| e^{-2ri\sin(\alpha)}\left(
\left(1+\frac{nt^3}{12}\right)\sum_{j=2}^{\infty}\left(\frac{t}{2}\right)^j\right)\right|d\alpha 
\leq \frac{\sqrt{e}}{2\pi }\int_{-\pi}^{\pi} 1.0003\cdot \frac{1}{4r^2(1-\frac{1}{2r})}d\alpha< \frac{0.413}{n}
\]
while the rightmost one gives
\[
\frac{\sqrt{e}}{2\pi }\int_{-\pi}^{\pi} \left| e^{-2ri\sin(\alpha)}\frac{nt^4}{24}\right|d\alpha=
\frac{\sqrt{e}}{2\pi }\cdot 2\pi \cdot \frac{n}{24n^2}=\frac{\sqrt{e}}{24n}.
\]
Therefore
\begin{equation}\label{E5}
E_5 \le \frac{0.413}{n}+ \frac{\sqrt{e}}{24n}.
\end{equation}

Let us now look at the integral $E_6$. We have
\begin{multline*}
E_6=\frac{\sqrt{e}}{2\pi }\int_{-\pi}^{\pi} e^{-2ri\sin(\alpha)}\left(\frac{nt^3}{12}+\frac{t}{2}\right)d\alpha
= \frac{\sqrt{e}}{2\pi }\int_{-\pi}^{\pi} e^{-2ri\sin(\alpha)}\left(\frac{e^{-3i\alpha}}{12r}+\frac{e^{-i\alpha}}{2r}\right)d\alpha\\
=\frac{\sqrt{e}}{2\pi r}\left(\frac{1}{12}\int_{-\pi}^{\pi} e^{-2ri\sin(\alpha)-3i\alpha}
d\alpha+\frac{1}{2}\int_{-\pi}^{\pi} e^{-2ri\sin(\alpha)-i\alpha}d\alpha\right)\\ 
=\frac{\sqrt{e}}{r}\left(\frac{1}{12}J_{-3}(2r)+\frac{1}{2}J_{-1}(2r)\right)=
\frac{\sqrt{e}}{2r}\left(-J_{1}(2r) - \frac{1}{6}J_3(2r)\right),
\end{multline*}
where $J_1$ and $J_3$ are J-Bessel functions.

Here we will apply Watson's treatise \cite{watson}, p. 205--206. 
Let $\nu\in\mathbb{Z}$. We assume that $\left(1+z\right)^{\nu-1/2}$ denotes the principal branch of the
multivalued function in the complex variable $z$.
Denote then
$$
P(x,v) = \frac{1}{\Gamma(v+1/2)}\int_0^{\infty}e^{-u}u^{v-1/2}p(u,x,v)du,\quad
p(u,x,v) = \frac{1}{2}\left(\left(1+\frac{iu}{2x}\right)^{v-1/2} + \left(1-\frac{iu}{2x}\right)^{v-1/2}\right),
$$
$$
Q(x,v) = \frac{1}{\Gamma(v+1/2)}\int_0^{\infty}e^{-u}u^{v-1/2}q(u,x,v)du,\quad
q(u,x,v) = \frac{1}{2i}\left(\left(1+\frac{iu}{2x}\right)^{v-1/2} - \left(1-\frac{iu}{2x}\right)^{v-1/2}\right),
$$
where $v\in\mathbb{Z}_{\ge 0}$ and $u,x\in\mathbb{R}$. 
Hence $p(u,x,v),q(u,x,v):\ \mathbb{R}\to\ \mathbb{R}$ and therefore we may apply the real domain Taylor's theorem 
$$
g(u) = \sum_{k=0}^{n} \frac{g^{(k)}(0)}{k!} u^k + \frac{g^{(n+1)}(s)}{(n+1)!} u^{n+1},\quad s\in [0,u],
$$
to the real valued functions $g(u)=p(u,x,v)$ and $g(u)=q(u,x,v)$, respectively. 

Watson writes the $v$th Bessel function as
\[
J_{v}(x )=\left(\frac{2}{\pi x}\right)^{1/2} \left( \cos\left( x - \frac{v\pi}{2} - \frac{\pi}{4} \right) P(x,v)
  - \sin\left( x - \frac{v\pi}{2} - \frac{\pi}{4} \right) Q(x,v) \right).
\]
Thus
\[
J_1(2r)=\left(\frac{1}{\pi r}\right)^{1/2}\left(\cos\left(2r-\frac{3\pi}{4}\right)P(2r,1) - \sin\left(2r-\frac{3\pi}{4}\right)Q(2r,1)\right),
\]
which we investigate first. 
First we note
$$
\frac{d}{dz}\left(1+az\right)^{v-1/2} = (v-1/2)a \left(1+az\right)^{v-3/2},\quad
\frac{d^2}{d^2z}\left(1+az\right)^{v-1/2} = (v-1/2)(v-3/2)a^2 \left(1+az\right)^{v-5/2},\quad$$
$$\frac{d^3}{d^3z}\left(1+az\right)^{v-1/2} = (v-1/2)(v-3/2)(v-5/2)a^3 \left(1+az\right)^{v-7/2}.
$$
Let $u\in[0,\infty]$. Set $a=i/4r$ and $z=u$. Then
\begin{equation*}
\begin{split}
\frac{d}{du}p(u,2r,1) & = \frac{d}{du} \frac{1}{2}\left( \left(1+\frac{iu}{4r}\right)^{1/2} + \left(1-\frac{iu}{4r}\right)^{1/2} \right)
              = \frac{1}{2}\frac{i}{4r}\cdot \frac{1}{2}\left( \left(1+\frac{iu}{4r}\right)^{-1/2} - \left(1-\frac{iu}{4r}\right)^{-1/2} \right), \\
\frac{d^2}{d^2u}p(u,2r,1) & = \frac{1}{2}\frac{-1}{2}\left(\frac{i}{4r}\right)^2
                            \cdot\frac{1}{2}\left( \left(1+\frac{iu}{4r}\right)^{-3/2} + \left(1-\frac{iu}{4r}\right)^{-3/2} \right),  \\
\frac{d^3}{d^3u}p(u,2r,1) & = \frac{1}{2}\frac{-1}{2}\frac{-3}{2} \left(\frac{i}{4r}\right)^3
                              \cdot\frac{1}{2}\left( \left(1+\frac{iu}{4r}\right)^{-5/2} - \left(1-\frac{iu}{4r}\right)^{-5/2} \right)
                           := \frac{3\Theta(u)}{2^9r^3}, \\																					
\end{split}
\end{equation*}
where $|\Theta(u)|\le 1$.  Thus
\begin{equation*}
\begin{split}
p(u,2r,1) & = p(0,2r,1) + \frac{p'(0,2r,1)}{1!}u + \frac{p''(0,2r,1)}{2!}u^2 + \frac{p'''(s,2r,1)}{3!}u^3 \\
          & = 1 + \frac{1}{2^7r^2}u^2 + \frac{\Theta_1(s)}{2^{10}r^3}u^3,\quad s\in [0,u],
\end{split}
\end{equation*}
where $|\Theta_1(s)|\le 1$ and $\Theta_1(u)\in\mathbb{R}$ (the same is valid for all $\Theta_i(u)$, $i=1,\ldots,6$). Further,
\begin{equation*}
\begin{split}
P(2r,1) & = \frac{1}{\Gamma(3/2)}\int_0^{\infty}e^{-u}u^{1/2} p(u,2r,1) du  \\
        & = \frac{1}{\Gamma(3/2)} \left( \Gamma(3/2) + \frac{1}{2^7r^2}\Gamma(7/2) + \frac{\Theta(u)}{2^{10}r^3} \Gamma(9/2) \right)
				  = 1 + \frac{15}{2^9r^2} + \frac{3\cdot 5\cdot 7}{2^{13}r^3}\,\Theta_1(s).
\end{split}
\end{equation*}

Next we look at the expression $Q(2r,1)$. Using a similar approach, we obtain
\begin{equation*}
\begin{split}
q(u,2r,1)  & = q(0,2r,1) + \frac{q'(0,2r,1)}{1!}u + \frac{q''(0,2r,1)}{2!}u^2 + \frac{q'''(s,2r,1)}{3!}u^3 \\
           & = \frac{1}{8r}u + \frac{\Theta_2(s)}{2^{10}r^3}u^3,\quad s\in [0,u],
\end{split}
\end{equation*}
where $|\Theta_2(s)|\le 1$. Thus
\begin{equation*}
\begin{split}
Q(2r,1) & = \frac{1}{\Gamma(3/2)}\int_0^{\infty} e^{-u}u^{1/2} q(u,2r,1) du \\
        & = \frac{1}{\Gamma(3/2)} \left( \frac{1}{8r} \Gamma(5/2) + \frac{\Theta_2(s)}{2^{10}r^3}\Gamma(9/2) \right)
				= \frac{3}{16r} + \frac{3\cdot 5\cdot 7}{2^{13}r^3}\,\Theta_2(s).
\end{split}
\end{equation*}

We may now move to 
\[
J_{3}(2r)=\left(\frac{1}{\pi r}\right)^{1/2} 
\left( \cos\left( 2r - \frac{7\pi}{4} \right) P(2r,3)  - \sin\left( 2r - \frac{7\pi}{4} \right) Q(2r,3) \right).
\]
Here
$$
\frac{d}{du}p(u,2r,3) = \frac{d}{du}\cdot\frac{1}{2}\left( \left(1+\frac{iu}{4r}\right)^{5/2} + \left(1-\frac{iu}{4r}\right)^{5/2} \right)
= \frac{5}{2}\frac{i}{4r}
\cdot\frac{1}{2}\left(\left(1+\frac{iu}{4r}\right)^{3/2} - \left(1-\frac{iu}{4r}\right)^{3/2} \right),
$$
$$
\frac{d^2}{d^2u}p(u,2r,3) = \frac{5}{2}\frac{3}{2}\left(\frac{i}{4r}\right)^2
\cdot\frac{1}{2}\left(\left(1+\frac{iu}{4r}\right)^{1/2} + \left(1-\frac{iu}{4r}\right)^{1/2} \right),
$$
$$
\frac{d^3}{d^3u}p(u,2r,3) = \frac{5}{2}\frac{3}{2}\frac{1}{2}\left(\frac{i}{4r}\right)^3
 \cdot\frac{1}{2}\left(\left(1+\frac{iu}{4r}\right)^{-1/2} - \left(1-\frac{iu}{4r}\right)^{-1/2} \right).
$$
Thus
\begin{equation*}
\begin{split}
p(u,2r,3) & = p(0,2r,3) + \frac{p'(0,2r,3)}{1!}u + \frac{p''(0,2r,3)}{2!}u^2 + \frac{p'''(s,2r,3)}{3!}u^3 \\
  & = 1 - \frac{15}{2^7r^2}u^2 + \frac{5\cdot\Theta_3(s)}{2^{10}r^3}u^3,\quad s\in [0,u],
\end{split}
\end{equation*}
where $|\Theta_3(s)|\le 1$. Here it was important to compute the third derivative in order to get negative exponents 
for the terms $1 \pm \frac{iu}{4r}$. This allows to receive a reasonable upper estimate for $\Theta_3(s)$.
Thereby
\begin{equation*}
\begin{split}
P(2r,3) & = \frac{1}{\Gamma(7/2)}\int_0^{\infty}e^{-u}u^{5/2}p(u,2r,3)du  \\
        & = \frac{1}{\Gamma(7/2)} \left( \Gamma(7/2) - \frac{15}{2^7r^2}\Gamma(11/2) + \frac{5\cdot\Theta_3(s)}{2^{10}r^3}\Gamma(13/2) \right) \\
				& = 1 - \frac{3^3\cdot 5\cdot 7}{2^9r^2} + \frac{3^2\cdot 5 \cdot 7\cdot 11}{2^{13}r^3}\,\Theta_3(s).
\end{split}
\end{equation*}
Now it suffices to look at the expression $Q(2r,3)$. Again, using a similar approach we get
\begin{equation*}
\begin{split}
q(u,2r,3) & = q(0,2r,3) + \frac{q'(0,2r,3)}{1!}u + \frac{q''(0,2r,3)}{2!}u^2  + \frac{q'''(s,2r,3)}{3!}u^3 \\
          & = \frac{5}{8r}u - \frac{5\cdot\Theta_4(s)}{2^{10}r^3}u^3,\quad s\in [0,u],
\end{split}
\end{equation*}where $|\Theta_4(s)|\le 1$. 
Consequently
\begin{equation*}
\begin{split}
Q(2r,3) & = \frac{1}{\Gamma(7/2)}\int_0^{\infty}e^{-u}u^{5/2} q(u,2r,3) du  \\
        & = \frac{1}{\Gamma(7/2)}\left( \frac{5}{8r}\Gamma(9/2) - \frac{5\cdot\Theta_4(s)}{2^{10}r^3}\Gamma(13/2) \right)  \\
			  & = \frac{5\cdot 7}{2^4r} - \frac{3^2\cdot 5 \cdot 7\cdot 11}{2^{13}r^3}\, \Theta_4(s).
\end{split}
\end{equation*}

Thus, we are led to estimate the integral

\[
H_1 = \frac{\sqrt{e}}{2\pi }\int_{-\pi}^{\pi} e^{-2ri\sin(\alpha)}d\alpha=\frac{\sqrt{e}}{2\pi }\cdot 2\pi J_{0}(2r)=\sqrt{e}J_{0}(2r).
\]
Using again \cite{watson}, p. 206, we have
\[
J_0(2r)=\left(\frac{1}{\pi r}\right)^{1/2}\left(\cos\left(2r-\frac{\pi}{4}\right)P(2r,0)-\sin\left(2r-\frac{\pi}{4}\right)Q(2r,0)\right).
\]
Here
\begin{equation*}
p(u,2r,0) = p(0,2r,0) + \frac{p'(0,2r,0)}{1!}u + \frac{p''(s,2r,0)}{2!}u^2 
          = 1 + \frac{3\Theta_5(s)}{2^7r^2}u^2,\quad s\in [0,u],
\end{equation*}
where $|\Theta_5(s)|\le 1$. Further,
\begin{equation*}
\begin{split}
P(2r,0) & = \frac{1}{\Gamma(1/2)}\int_0^{\infty}e^{-u}u^{-1/2}p(u,2r,0)du  \\
        & = \frac{1}{\Gamma(1/2)} \left( \Gamma(1/2) + \frac{3\Theta_5(s)}{2^7r^2}\Gamma(5/2) \right)
				= 1 + \frac{3^2\Theta_5(s)}{2^9r^2}.
\end{split}
\end{equation*}

For the last we consider $Q(2r,0)$. Now
\begin{equation*}
q(u,2r,0)  = q(0,2r,0) + \frac{q'(0,2r,0)}{1!}u + \frac{q''(s,2r,0)}{2!}u^2  
            = -\frac{1}{8r}u  + \frac{3\Theta_6(s)}{2^7r^2}u^2,\quad s\in [0,u],
\end{equation*}
where $|\Theta_6(s)|\le 1$. Further,
\begin{equation*}
\begin{split}Q(2r,0) & = \frac{1}{\Gamma(1/2)}\int_0^{\infty} e^{-u}u^{-1/2} q(u,2r,0)du  \\
                     & = \frac{1}{\Gamma(1/2)} \left( -\frac{1}{8r}\Gamma(3/2) +\frac{3\Theta_6(s)}{2^7r^2}\Gamma(5/2)\right)
				               =  -\frac{1}{16r} + \frac{3^2\Theta_6(s)}{2^9r^2}.
\end{split}
\end{equation*}
Hence,
\[
H_1 =\sqrt{e}J_{0}(2r)=\sqrt{\frac{e}{\pi r}}\left(\cos\left(2r-\frac{\pi}{4}\right) 
+ \frac{1}{16r}\cdot\sin\left(2r-\frac{\pi}{4}\right)+\frac{3^2\Theta_7}{2^8r^2}\right),\quad |\Theta_7|\le 1.
\]
Now we want to put all the Bessel function terms together. Observe first that

\begin{equation*}
\begin{split}
I_4 - E_5 
    & \le H_1+E_6 = \sqrt{e}J_{0}(2r) - \frac{\sqrt{e}}{2r}J_1(2r)-\frac{\sqrt{e}}{12r}J_3(2r) \\
    & = \sqrt{\frac{e}{\pi}}\left( \frac{1}{r^{1/2}} \left(\cos\left(2r-\frac{\pi}{4}\right)\left(1 + \frac{3^2\Theta_5(s)}{2^9r^2}\right)
		   - \sin\left(2r-\frac{\pi}{4}\right)\left(-\frac{1}{16r} + \frac{3^2\Theta_6(s)}{2^9r^2}\right) \right) \right.\\
		& - \frac{1}{2r^{3/2}} 
		    \left(\cos\left(2r-\frac{3\pi}{4}\right) \left(1 + \frac{15}{2^9r^2} + \frac{3\cdot 5\cdot 7}{2^{13}r^3}\,\Theta_1(s) \right)
		        - \sin\left(2r-\frac{3\pi}{4}\right)\left(\frac{3}{16r} + \frac{3\cdot 5\cdot 7}{2^{13}r^3}\,\Theta_2(s)\right) \right)  \\
    & - \frac{1}{12r^{3/2}} 
		    \left( \cos\left( 2r - \frac{7\pi}{4} \right) \left(1 - \frac{3^3\cdot 5\cdot 7}{2^9r^2} 
				                     + \frac{3^2\cdot 5 \cdot 7\cdot 11}{2^{13}r^3}\,\Theta_3(s) \right) \right.\\
		&  \left.  - \sin\left( 2r - \frac{7\pi}{4} \right) \left(\frac{5\cdot 7}{2^4r} 
						            - \frac{3^2\cdot 5 \cdot 7\cdot 11}{2^{13}r^3}\, \Theta_4(s) \right) \right).														
\end{split}
\end{equation*}
The main term, $\sqrt{\frac{e}{\pi}}\frac{1}{r^{1/2}} \cos\left(2r-\frac{\pi}{4}\right)$,
in $H_1$ yields the main term in the whole $I_4-E_5$.
Let us now concentrate on those terms which have the term $r^{-3/2}$. Putting all of these together, we obtain
\begin{equation*}
\begin{split}
&  \sqrt{\frac{e}{\pi}}\,\frac{1}{r^{3/2}}
   \left(\frac{1}{16}\sin\left(2r-\frac{\pi}{4}\right) - \frac{1}{2}\cos\left(2r-\frac{3\pi}{4}\right) 
	 - \frac{1}{12}\cos\left(2r-\frac{7\pi}{4}\right) \right)\\
& = \sqrt{\frac{e}{\pi}}\,\frac{1}{r^{3/2}}
    \left(\frac{1}{16}\sin\left(2r-\frac{\pi}{4}\right) - \frac{1}{2}\sin\left(2r-\frac{\pi}{4}\right) 
		+ \frac{1}{12}\sin\left(2r-\frac{\pi}{4} \right) \right)\\
& = -\frac{17}{48} \sqrt{\frac{e}{\pi}}\cdot\frac{1}{r^{3/2}}\sin\left(2r-\frac{\pi}{4}\right).
\end{split}
\end{equation*}
Let us now move to the rest of the terms in $I_4-E_5$. We use absolute values to bound the terms. Their contribution is at most
\begin{equation}\label{}
\begin{split}
    &  \sqrt{\frac{e}{\pi}}\ \Biggl| \frac{1}{r^{1/2}} \left(\cos\left(2r-\frac{\pi}{4}\right)\frac{3^2\Theta_5(s)}{2^9r^2}
		     - \sin\left(2r-\frac{\pi}{4}\right) \frac{3^2\Theta_6(s)}{2^9r^2} \right)  \\
		& - \frac{1}{2r^{3/2}} 
		    \left(\cos\left(2r-\frac{3\pi}{4}\right) \left(\frac{15}{2^9r^2} + \frac{3\cdot 5\cdot 7}{2^{13}r^3}\,\Theta_1(s) \right)
		        - \sin\left(2r-\frac{3\pi}{4}\right)\left(\frac{3}{16r} + \frac{3\cdot 5\cdot 7}{2^{13}r^3}\,\Theta_2(s)\right) \right)  \\
    & - \frac{1}{12r^{3/2}} 
		    \biggl( \cos\left( 2r - \frac{7\pi}{4} \right) \left( - \frac{3^3\cdot 5\cdot 7}{2^9r^2} 
				                     + \frac{3^2\cdot 5 \cdot 7\cdot 11}{2^{13}r^3}\,\Theta_3(s) \right)  \\
		&		    - \sin\left( 2r - \frac{7\pi}{4} \right) \left(\frac{5\cdot 7}{2^4r} 
						            - \frac{3^2\cdot 5 \cdot 7\cdot 11}{2^{13}r^3}\, \Theta_4(s) \right) \biggr) \Biggr|,   
\end{split}
\end{equation}
which is estimated by												
\begin{equation}\label{E7}
\begin{split}												
\le & \sqrt{\frac{e}{\pi}}\, \Biggl( \frac{1}{r^{1/2}} \left( \frac{3^2}{2^9r^2} + \frac{3^2}{2^9r^2} \right) 
		  + \frac{1}{2r^{3/2}} 
			\left( \frac{15}{2^9r^2} + \frac{3\cdot 5\cdot 7}{2^{13}r^3} + \frac{3}{16r} + \frac{3\cdot 5\cdot 7}{2^{13}r^3} \right) \\
  + & \frac{1}{12r^{3/2}} 
		    \left( \frac{3^3\cdot 5\cdot 7}{2^9r^2} + \frac{3^2\cdot 5 \cdot 7\cdot 11}{2^{13}r^3} 
						   + \frac{5\cdot 7}{2^4r} + \frac{3^2\cdot 5 \cdot 7\cdot 11}{2^{13}r^3} \right) \Biggr) \\ 																		
  = & \sqrt{\frac{e}{\pi}}\, \Biggl( \frac{239}{768n^{5/4}} + \frac{345}{2048n^{7/4}} + \frac{1365}{16384n^{9/4}} \Biggr) =: E_7.
\end{split}
\end{equation}

\subsection{Finishing the proof: putting all the bounds together}

The main term in $H_1$ yields the main term. For the error term, it suffices to put all errors together. 
By \eqref{E2}, \eqref{E3}, \eqref{E4}, \eqref{E5}, \eqref{E7} and by using the bound $\frac{1}{\sqrt{n}-1/2}<\frac{1.002}{\sqrt{n}}$, 
which is valid for $n\geq 10^5$, we get
\begin{equation*}\label{}
\begin{split}	
& \left|L_n(1)-\sqrt{\frac{e}{\pi}}\cdot \frac{1}{n^{1/4}}\cos\left(2\sqrt{n}-\frac{\pi}{4}\right)
+\frac{17}{48} \sqrt{\frac{e}{\pi}}\cdot\frac{1}{n^{3/4}}\sin\left(2\sqrt{n}-\frac{\pi}{4}\right)\right| \\ 
& \leq E_2+E_3+E_4+E_5+E_7 \\ 
& < \frac{0.00605\cdot e^{3/2}}{n^{3/2}(\sqrt{n}-1/2)}
   +\frac{ e^{3/2}}{80n(\sqrt{n}-1/2)}+\frac{0.007\cdot\sqrt{e}}{n^{1/2}(\sqrt{n}-1/2)}+\frac{0.413}{n}+\frac{\sqrt{e}}{24n} \\
&  +\sqrt{\frac{e}{\pi}}\, \Biggl( \frac{239}{768n^{5/4}} + \frac{345}{2048n^{7/4}} + \frac{1365}{16384n^{9/4}} \Biggr) \\ 
& < \frac{0.00605\cdot e^{3/2}\cdot 1.002}{10^5n}
   +\frac{ e^{3/2}\cdot 1.002}{80n\sqrt{10^5}}+\frac{0.007\cdot\sqrt{e}\cdot 1.002}{n}+\frac{0.413}{n}
   +\frac{\sqrt{e}}{24n} \\ 
& + \sqrt{\frac{e}{\pi}}\, \Biggl( \frac{239}{768n(10^5)^{1/4}} + \frac{345}{2048n (10^5)^{3/4}} + \frac{1365}{16384n(10^5)^{5/4}} \Biggr)
  < \frac{0.510}{n}.
\end{split}
\end{equation*}

This ends the proof of Theorem \ref{laguerre2}.

Finally, noticing that
\[
|L_n(1)|\leq \left|\sqrt{\frac{e}{\pi}}\cdot \frac{\cos (2\sqrt{n}-\frac{\pi}{4})}{n^{1/4}}\right|
+\left|\frac{17}{48}\sqrt{\frac{e}{\pi}}\cdot\frac{\sin(2\sqrt{n}-\frac{\pi}{4})}{n^{3/4}}\right|+\frac{0.510}{n}<1
\]
for $n\geq 4$, and numerically we can check that $|L_n(1)|\leq 1$ for $n\in \{1,2,3\}$, we obtain Corollary \ref{laguerre3}.


\begin{thebibliography}{}

\bibitem{BBC2008} Borwein D., Borwein J.M. and Crandall R.E., Effective Laguerre Asymptotics. SIAM Journal on Numerical Analysis 
46 (2008), 3285-3312. 

\bibitem{AMLOTA2023}  Ernvall-Hyt\"onen A-M., Matala-aho T. and Sepp\"al\"a L.,
Euler's factorial series, Hardy integral, and continued fractions. J. Number Theory 244 (2023), 224-250.

\bibitem{Fejer1909} Fej\'er L., Asymptotikus \'ert\'ekek meghat\'aroz\'as\'arol. Mathematikai \'es Term\'eszettudom\'anyi
\'Ertesit\"o, 27 (1909), 1-33.

\bibitem{GradshteynRyzhik}  Gradshteyn I.S. and Ryzhik I.M., Table of Integrals, Series, and Products.
Academic Press, 7th ed., 2007.

\bibitem{KoomanTijdeman1990} Kooman R.J., Tijdeman R., Convergence properties of linear recurrence sequences. 
Nieuw Arch. Wiskd., IV. Ser. 8, No. 1, (1990), 13-25.

\bibitem{korolev} Korolev, M.A., Gram's Law in the Theory of the Riemann Zeta-Function, Part 1.
Proc. Steklov Inst. Math. 292 (Suppl. 2) (2016), 1-146.

\bibitem{TAWA2018} Matala-aho T., Zudilin W., Euler's factorial series and global relations. J. Number Theory 186 (2018), 202-210.

\bibitem{sage} SageMath, the Sage Mathematics Software System (Version 8.9). The Sage Developers, 2019, https://www.sagemath.org.	

\bibitem{Szego1975} Szeg\"o G., Orthogonal Polynomials. Amer. Math. Soc. Colloq. Publ., 4th ed., 23, American
Mathematical Society, Providence, RI, 1975.

\bibitem{Tricomi1949} Tricomi F., Sul comportamento asintotico dei polinomi di Laguerre. Annali di Matematica Pura ed Applicata, 28 (1949), 263-289.

\bibitem{watson} Watson, G. N., A treatise on the theory of Bessel functions. Cambridge University Press, 1966.

\end{thebibliography}
\end{document}